\documentstyle{article}
\begin{document}
\newtheorem{tm}{Theorem}
\newtheorem{lm}{Lemma}
\newtheorem{df}{Definition}
\newtheorem{rem}{Remark}
\newtheorem{prop}{Proposition}
\newtheorem{cor}{Corollary}
\newtheorem{ex}{Example}
\newcommand {\Oe} {\Omega_\varepsilon}
\newcommand {\Oet} {\Omega_\varepsilon^T}
\newcommand {\f} {\varphi}
\newcommand   {\ep}  {\varepsilon}
\newcommand {\fe} {\varphi_\varepsilon}

\begin{center}
{\large MATHEMATICAL MODELING OF AN ARRAY  OF NUCLEAR WASTE CONTAINERS}\\
\ \\
{\bf Alain Bourgeat\footnote{MCS-ISTIL, Universit\'{e} Lyon1, B\^at. ISTIL, 43 Bd. du 11
 novembre, 69622 Villeurbanne Cedex, France} ,
 Olivier Gipouloux\footnote{Facult\'{e} de sciences, Universit\'{e} de St-Etienne, 23 Rue
 Dr.Paul Michelain, St-Etienne Cedex 2, France}\footnote{Laboratoire de M\'ecanique et
 d'Acoustique, UPR  7051,31 Chemin Joseph Aiguier, 13402 Marseille cedex 20, France}
 and Eduard Maru\v{s}i\'{c}-Paloka\footnote{ Department of Mathematics, University
 of Zagreb, Bijeni\v{c}ka 30, 10000 Zagreb, Croatia}}
\ \\
\end{center}\section{Introduction}

The goal of this paper is to give a mathematical model describing  the global behavior
of
 an underground waste repository, once the containers start to leak.
The purpose of such a global model is to be used for the full field simulations used
in
 safety assessements.
The physical situation can be described as an array made of high number of leaking
modules inside a thin low permeable layer (e.g. clay), included between two bigger layers
 with higher permeability (e.g. limestone or marl). The pollutant is transported both by
the convection produced by the water flowing slowly (creeping flow)through
the rocks and by the diffusion coming from the dilution in the water. The leaking last
the all period of time $]0,t_m[$, that is small compared to the millions of years over
which convection and diffusion are active. In a real repository there is a pressure drop
 producing the flow crossing a large number of disposal modules where each module includes
 several containers. Herein, for simplicity, the repository consists of a set of modules
lying on a hypersurface $\Sigma$ and we represent the leaking of a disposal module by a localized
density source inside the domain or by a hole in the domain with  a given flux on its boundary.
Moreover, without lost of generality we assume the convection velocity field
to be given . According to the test case \cite{CPLEX}, the typical size of a  module is a
hundred of meters for the width, a kilometer for the length and five meters for the height.
The distance between two modules is also of order 100 meters and the low permeable layer
(the clay layer), in which  the repository is embedded, has respectively a height and a
length  of order 150 and 3000 meters.
Since there is a large number of modules, each of them with a small size compared to the
layers size, see figure 1, a direct numerical simulations of the full field, based on a
{\em microscopic} model taking in account all the detail, is unrealistic. Considering
 the ratio
 between the width of a single module $l$ and the layer length $L$, which is of order $1/30 $,
 as a small parameter, $\varepsilon $, in the {\em microscopic} model, then
 the modules, have a height of order $\varepsilon ^2$, and are now
inmbedded in a layer of thickness $\varepsilon$.
 The study of the renormalized model behavior, as $\varepsilon $ tends to $0$
, by means of the homogenization
method and boundary layers, gives an asymptotic model which could be used as a repository
global model for numerical simulations.\\
\ \\
\ \\
\ \\
\ \\
\ \\
\ \\
\ \\
\begin{picture}(60,40)(-70,130)
\setlength{\unitlength}{1 mm}
\put(20,80){\line(1,0){60}}
\put(20,20){\line(1,0){60}}
\put(20,20){\line(0,1){60}}
\put(80,20){\line(0,1){60}}
\put(20,48){\line(1,0){60}}
\put(20,52){\line(1,0){60}}
\put(21,49.8){\line(1,0){2}}

\put(21,49.8){\line(1,0){2}}
\put(21,49.8){\line(0,1){0.4}}
\put(23,49.8){\line(0,1){0.4}}

\put(25,49.8){\line(1,0){2}}
\put(25,49.8){\line(0,1){0.4}}
\put(27,49.8){\line(0,1){0.4}}

\put(29,49.8){\line(1,0){2}}
\put(29,49.8){\line(0,1){0.4}}
\put(31,49.8){\line(0,1){0.4}}

\put(33,49.8){\line(1,0){2}}
\put(33,49.8){\line(0,1){0.4}}
\put(35,49.8){\line(0,1){0.4}}

\put(37,49.8){\line(1,0){2}}
\put(37,49.8){\line(0,1){0.4}}
\put(39,49.8){\line(0,1){0.4}}

\put(41,49.8){\line(1,0){2}}
\put(41,49.8){\line(0,1){0.4}}
\put(43,49.8){\line(0,1){0.4}}

\put(45,49.8){\line(1,0){2}}
\put(45,49.8){\line(0,1){0.4}}
\put(47,49.8){\line(0,1){0.4}}

\put(49,49.8){\line(1,0){2}}
\put(49,49.8){\line(0,1){0.4}}
\put(51,49.8){\line(0,1){0.4}}

\put(53,49.8){\line(1,0){2}}
\put(53,49.8){\line(0,1){0.4}}
\put(55,49.8){\line(0,1){0.4}}

\put(57,49.8){\line(1,0){2}}
\put(57,49.8){\line(0,1){0.4}}
\put(59,49.8){\line(0,1){0.4}}

\put(61,49.8){\line(1,0){2}}
\put(61,49.8){\line(0,1){0.4}}
\put(63,49.8){\line(0,1){0.4}}

\put(65,49.8){\line(1,0){2}}
\put(65,49.8){\line(0,1){0.4}}
\put(67,49.8){\line(0,1){0.4}}

\put(69,49.8){\line(1,0){2}}
\put(69,49.8){\line(0,1){0.4}}
\put(71,49.8){\line(0,1){0.4}}

\put(73,49.8){\line(1,0){2}}
\put(73,49.8){\line(0,1){0.4}}
\put(75,49.8){\line(0,1){0.4}}

\put(77,49.8){\line(1,0){2}}
\put(77,49.8){\line(0,1){0.4}}
\put(79,49.8){\line(0,1){0.4}}

\put(21,50.2){\line(1,0){2}}

\put(25,50.2){\line(1,0){2}}

\put(29,50.2){\line(1,0){2}}

\put(33,50.2){\line(1,0){2}}
\put(33,45){\line(0,1){5}}
\put(35,45){\line(0,1){5}}

\put(34,47){\vector(-1,0){1}}
\put(34,47){\vector(1,0){1}}
\put(30,42){$\textrm{{\footnotesize 100 m }}$}

\put(37,50.2){\line(1,0){2}}

\put(41,50.2){\line(1,0){2}}

\put(45,50.2){\line(1,0){2}}

\put(49,50.2){\line(1,0){2}}

\put(53,50.2){\line(1,0){2}}

\put(57,50.2){\line(1,0){2}}

\put(61,50.2){\line(1,0){2}}

\put(65,50.2){\line(1,0){2}}

\put(69,50.2){\line(1,0){2}}

\put(73,50.2){\line(1,0){2}}

\put(77,50.2){\line(1,0){2}}

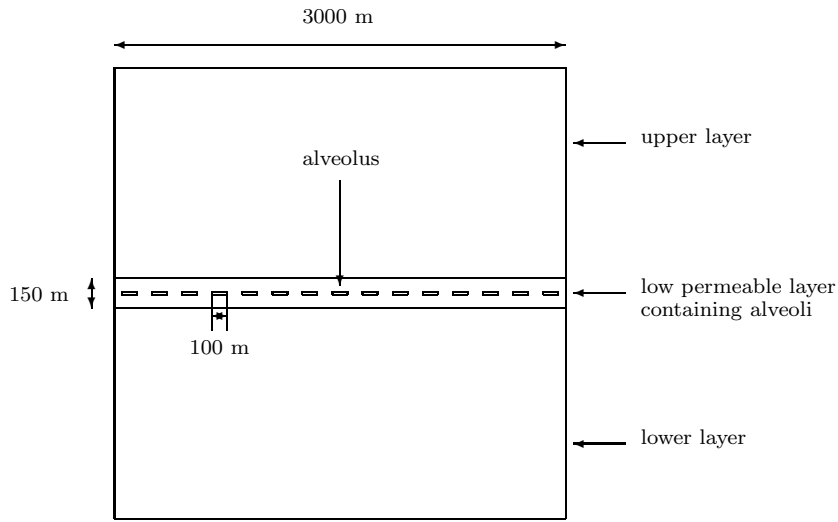
\put(20,12){$\textrm{{\footnotesize Figure 1: Three layers of soil containing alveoli}}$}
\put(88,50){\vector(-1,0){7}}
\put(90,30){$\textrm{{\footnotesize lower layer}}$}
\put(88,30){\vector(-1,0){7}}
\put(90,50){$\textrm{{\footnotesize low permeable layer }}$}
\put(90,47){$\textrm{{\footnotesize containing alveoli }}$}

\put(88,70){\vector(-1,0){7}}
\put(90,70){$\textrm{{\footnotesize upper layer}}$}
\put(50,65){\vector(0,-1){14}}
\put(45,67){$\textrm{{\footnotesize alveolus }}$}
\put(17,50){\vector(0,1){2}}
\put(17,50){\vector(0,-1){2}}
\put(6,49){$\textrm{{\footnotesize 150 m }}$}
\put(50,83){\vector(1,0){30}}
\put(50,83){\vector(-1,-0){30}}
\put(45,86){$\textrm{{\footnotesize 3000 m }}$}

\end{picture} \\ 
\ \\
\ \\
\ \\
\ \\
\ \\
\ \\  
\ \\
\ \\
\ \\
 We use methods similar to those applied to the fluid flow through a sieve in \cite{CC}, \cite{SP} or \cite{BGM}. Similar stationary problem with zero source term (i.e. $\Phi = 0 $) was treated in \cite{DV}.

\section{Setting the problem}
\subsection{ The Geometry}
Let $\Omega \subset {\bf R}^n $ be a bounded domain. Let
$ \Sigma \subset {\bf R}^{n-1}\,\times\,\{0\} $ be such that
$\Sigma \subset\Omega $. Let $A\subset {\bf R}^{n-1}$ be a periodic set
obtained by periodic repetition of bounded closed set $\,M\subset\subset\,]-1/2,1/2 [^{n-1}\;.$ More precisely
$$A=\bigcup_{\alpha \in {\bf Z}^{n-1}} M_\alpha \;\;\;,$$
where $M_\alpha=\alpha + M $. For small parameter $\varepsilon \ll 1 $ and
$\beta > 1 $ we define
$B_\varepsilon = (\,\varepsilon \;A\,\cap \Sigma\,)\,\times\,]\,
-\varepsilon^\beta , \varepsilon^\beta \,[  \,$ (in situation described in the introduction $\beta =2$). We denote by $J(\ep) =
\{\alpha \in {\bf Z}^{n-1}\;\;;\;\ep \;M_\alpha \cap \Sigma \neq \emptyset\;\}$
 and by $\Gamma_\alpha^\ep = \partial (\ep\,M_\alpha \,\times \,]-\ep^\beta ,
\ep^\beta [\;)\;$ . Finally $\Oe = \Omega \,\backslash\,B_\varepsilon $,
$\Oet = (\Omega \,\backslash\,B_\varepsilon ) \,\times \,]0,T[\;$, $\Omega^T =
\Omega \,\times \,]0,T[\;$ and
$\Gamma_\ep = \partial B_\varepsilon $, $\,\Gamma_\ep^T= \Gamma_\ep\,
\times\,]0,T[ $ .
\subsection{The Equations}
Let $\Phi\in L^\infty ([0,T])$ be the function describing the time behaviour of an alveolus. As mentioned before it has a compact support $[0,t_m] \subset ]0,T[ $. Let $\lambda =\frac{\log 2}{\tau}> 0 $, with $\tau $  being the half life of the radioactive element, and let $\f_0\in H^1(\Oe)$ the initial concentration of the radioactive material in the soil (typically equal to zero). The diffusion is described by
${\bf A}\in L^\infty ({\bf R};{\bf R}^{n\times n}) $ a positive definite matrix function. Since layers of soil involved in our model have different properties, we assume that 
$${\bf A}(y_n )=\left\{
\begin{array}{l}

{\bf A}^1\;\;\;,\;\;\mbox{for}\;\;|y_n|<h\\
{\bf A}^2\;\;\;,\;\;\mbox{for}\;\;|y_n|>h
\end{array}\right.$$
Now we write our diffusion matrix in the form ${\bf A}^\ep(x_n )={\bf A}(\frac{x_n}{\ep} )$. In the above described situation the low permeable layer has a hight of 150 m meaning that, in that case, $h=3/2 $.
We have the same situation with the convection velocity ${\bf v}\in C([0,T];H^1 ({\bf R}\times \Omega )^n)$. The dependence on $y_n$ is similar as in the case of diffusion matrix:
$${\bf v}(x,y_n,t )=\left\{
\begin{array}{l}
{\bf v}^1 (x,t)\;\;\;,\;\;\mbox{for}\;\;|y_n|<h\\
{\bf v}^2 (x,t)\;\;\;,\;\;\mbox{for}\;\;|y_n|>h
\end{array}\right.$$
For simplicity, we assume that the last component $v_n $ does not depend on $y_n$. Next, we suppose that $\mbox{div}_x{\bf v} =0 $, in order to have the divergence free convection velocity. Finally we pose ${\bf v}^\ep (x,t)={\bf v}(x,\frac{x_n}{\ep} ,t) $. At last, we define the porosity of the medium as
$${\bf \omega}(y_n )=\left\{
\begin{array}{l}
{\bf \omega}^1\;\;\;,\;\;\mbox{for}\;\;|y_n|<h\\
{\bf \omega}^2\;\;\;,\;\;\mbox{for}\;\;|y_n|>h
\end{array}\right.$$
and we put
$$\omega^\ep (x_n) = \omega (x_n /\ep )\;\;\;.$$
\\
The process is governed by the following convection-diffusion type equation:
\begin{eqnarray}
&&\omega^\ep\;\frac{\partial \fe}{\partial t} - \mbox{div}\,({\bf A}^\ep\nabla \fe ) +
({\bf v}^\ep\,\cdot\,\nabla\,)\fe+\lambda\,\omega^\ep\,\fe =0\;\;\;\mbox{in}\;\;\Oet
\label{prva,1} \\
&&\fe (0,x)=\f_0 (x)\;\;x\in \Oe \label{pu,1}\\
&&{\bf n}\,\cdot\, ({\bf A}^\ep\nabla \fe - {\bf v}^\ep\;\fe ) =\Phi (t)\;\;\mbox{on}\;\Gamma_\varepsilon^T \label{source,1}
\end{eqnarray}
We also need to impose some boundary condition on the exterior boundary
$S=\partial \Omega $.
Let $S=S_1 \cup S_2 $, where $S_i$ are disjoint and connected parts of $S$.
We impose
\begin{eqnarray}
&&\fe=0\;\;\;\;\mbox{on}\;\;S_1\\
&&{\bf n}\,\cdot\,{\bf A}^\ep\nabla \fe -({\bf v}^\ep\cdot{\bf n})\,\fe=0\;\;
\;\mbox{on}\;\;S_2\;\;.
\label{zadnja,1}
\end{eqnarray}
\section{ A priori estimates}
The main result of this section is:
\begin{prop}
Let $\fe $ be a unique solution of (\ref{prva,1})-(\ref{zadnja,1}). Then
there exists a constant $C>0 $  independent of $\ep$
such that
\begin{eqnarray}
&&|\fe |_{L^\infty (\Oet )} \leq C\label{mak,1}\\
&&|\fe |_{L^2 (0,T;H^1 (\Oe ))} \leq C\label{hi,1}
\end{eqnarray}\label{aprior,1}
\end{prop}
{\bf Proof.} The estimate (\ref{mak,1}) is the consequence of the maximum
principle. To prove (\ref{hi,1}) we use $\fe $ as the test function in
(\ref{prva,1})-(\ref{zadnja,1}). We obtain
\begin{eqnarray*}
&&\left|\frac{\sqrt{\omega^\ep} }{2}\;\fe (\,\cdot\, ,T ) \right|^2_{L^2 (\Oe)}+
({\bf A}^\ep \nabla\fe | \nabla\fe )_{L^2 (\Oet)}+
\lambda\;|\sqrt{\omega^\ep} \,\fe |^2_{L^2 (\Oet )}=\\
&&=
\int_0^T\;\Phi\;\sum_{\alpha\in J(\ep )}\,\int_{\Gamma_\alpha^\ep} \fe +
 \left|\;\frac{\sqrt{\omega^\ep}}{2}\;\f_0 \right|^2_{L^2 (\Oe )}\leq\\
&&\leq C\;(1+|\fe |_{H^1 (\Oet )})\;\;.\;\;\Box
\end{eqnarray*}
\section{Weak convergence}
Our solution $\fe $ is defined on variable domain $\Oet $. To use the weak convergence methods,
we extend it to  whole domain $\Omega^T $ preserving the estimates (\ref{mak,1}), (\ref{hi,1}).
In the sequel we assume that $\fe $
is extended using the results from \cite{DP} and we denote that extension by the same symbol.
Due to the proposition \ref{aprior,1}, we can conclude that there exists some
$\f \in L^2 (0,T;H^1(\Omega ))\cap L^\infty (\Omega^T)$
such that (up to a subsequence)
\begin{eqnarray}
&\fe \rightharpoonup \f\;\;\;\;&\;\;\;\mbox{weak* in}\;\;\;  L^\infty (\Omega^T)\\
&\nabla\fe \rightharpoonup \nabla\f \;\;\;&\;\;\;\mbox{weakly in}\;\;\;L^2 (\Omega^T)\;\;\;.
\end{eqnarray}
The main goal of this section is to identify the limit $\f$ .
We prove that:
\begin{tm}
The limit function $\f$ is the unique solution of the problem
\begin{eqnarray}
&&\omega^2\;\frac{\partial \f}{\partial t} - \mbox{div}\,({\bf A}^2\nabla \f ) +
({\bf v}^2\,\cdot\,\nabla\,)\f+\lambda\,\omega^2\,\f =0\;\;\;\mbox{in}\;\;
\tilde{\Omega}^T=(\Omega \backslash \Sigma )\,\times\,]0,T[
\label{prvak,1} \\
&&\f (x,0)=\f_0 (x)\;\;x\in \tilde{\Omega }=\Omega \backslash \Sigma \label{pocuv,1}\\
&&\f=0\;\;\;\;\mbox{on}\;\;S_1\label{rubuv,1}\\
&&{\bf n}\,\cdot\,{\bf A}^2\nabla \f -({\bf v}^2\cdot{\bf n})\,\f=0\;\;
\;\mbox{on}\;\;S_2\;\;\label{rubuvv,1}\\
&&[\f ] = 0\;\;\;,\;\;\;\left[{\bf e}_n\,\cdot\,{\bf A}^2\nabla \f
-({\bf v}^2\cdot{\bf e}_n)\,\f \right] =2\;\Phi\;|M |\;\;\mbox{on}\;\;\Sigma\;\;,
\label{zadnjak,1}
\end{eqnarray}
where $[w](x')=w(x',0+)-w(x',0-)\,,\,x'=(x_1,\ldots ,x_{n-1})$ denotes the jump over $\Sigma $ and $|M | $ denotes 
the area of $M $. 
\end{tm}
{\bf Proof.} Let $\psi \in C^\infty ([0,T];C^\infty_0 (\Omega )) $ 
be such that  $\psi (\;\cdot \;, T)= 0$ . Using $\psi $ as the test function in
(\ref{prva,1})-(\ref{zadnja,1}) we get
\begin{eqnarray*}
&&0=-\int_{\Oet }\omega^\ep\; \fe\;\frac{\partial\psi}{\partial t} - \int_{\Oe} \f_0\; 
\psi (\;\cdot \;,0) + \int_{\Oet} {\bf A}^\ep \nabla \fe \;\nabla\psi+\\
&&+\int_{\Oet} \omega^\ep\,\lambda\;\fe\;\psi+ \int_0^T
\Phi\; \sum_{i\in J(\ep )} \int_{\Gamma^\ep_i } \psi\;\;.
\end{eqnarray*}
Passage to the limit for the first four integrals is straightforward. 
For the last integral we have
$$\int_{\Gamma^\ep_i} \psi (x,t)\, d\Gamma^\ep_i = (\psi (x_i^\ep ,t) +O(\ep )\;) \;
|\Gamma^\ep_i |=\psi (x_i^\ep ,t) \;
2\;|M|\;\ep^{n-1} +O(\ep^{n+\beta -2})\;\;,$$
where $x_i^\ep =(\;(x_i^\ep )',0) $ is an arbitrary point from $\ep\,M_\alpha \,\times \{ 0\} $. But then
$$\sum_{i\in J(\ep )} \int_{\Gamma^\ep_i } \psi (x,t)\;d\Gamma^\ep_i \to 
 2\;|M|\;\int_\Sigma \psi (x',0,t)\;dx'\;\;,$$
where $x'=(x_1,\ldots ,x_{n-1})\in {\bf R}^{n-1} $ . $\;\;\Box$
\begin{rem}
In fact we did not use the periodicity of distribution of alveoli. The same proof holds in case
if each alveolus is randomly placed in a mash of an $\ep $ net. The alveoli do not even need to have the same shape, only the areas of their surfaces need to be equal.
\end{rem}
\section{Asymptotic expansion}
The above weak limit describes the global long time behaviour of the process in case when the flux $\Phi $ is not too large . However if we need more accurate information on the behaviour in the near field (i.e. in vicinity of $\Sigma) $, more precise asymptotics is needed.\\
To avoid cumbersome computations we simplify the geometry by assuming that $\Omega = ]-L/2 , L/2 [^n $. We denote then
\begin{eqnarray*}
&&\Sigma= ]-L/2 , L/2 [^{n-1} \times \{0\}\\
 &&S^+= ]-L/2 , L/2 [^{n-1} \times \{L/2\}\\
&&S^-= ]-L/2 , L/2 [^{n-1} \times \{-L/2\}\;\; .
\end{eqnarray*}
We also suppose that the alveoli are rectangular (which is true in real-world situation). More precisely, we take
$$M=\,\prod_{i=1}^{n-1} \,]-m_i , m_i\,[\;\;\;,\;\;1/2>m_i> 0\;\;.$$
We impose the Dirichlet condition on the bottom and Neumann condition on the top of the domain:

\begin{eqnarray}
&&\fe=0\;\;\;\;\mbox{on}\;\;S^-\\
&&{\bf n}\,\cdot\,{\bf A}^\ep\nabla \fe -({\bf v}^\ep\cdot{\bf n})\,\fe=0\;\;
\;\mbox{on}\;\;S^+\;\;.
\label{zadnjas,1}
\end{eqnarray}
On the lateral boundary we impose the periodicity condition in order to avoid the lateral boundary layer. More precisely, we impose
\begin{equation}
\fe\;\;\;\;\mbox{is}\;\;L-\mbox{periodic in}\;\;x'=(x_1,\ldots ,x_{n-1}) \label{per,1}
\end{equation}
We also assume that $L/\ep \in {\bf N} $. For the sake of compatibility, we suppose that given data $v^1 ,v^2$ and $\varphi_0$ are $L$-periodic in $x'$.\\
\ \\
We expect some fast changes of solution in vicinity of containers. Therefore, in that region, we introduce the
fast variable $y=x/\ep $ to describe that behaviour. Far from the sources we expect $\fe $,  the solution
of (\ref{prva,1})-(\ref{source,1}), (\ref{zadnjas,1}) and (\ref{per,1}),  to behave almost like our weak limit
$\f $, which, in this case, satisfies (\ref{prvak,1}), (\ref{pocuv,1}) , (\ref{zadnjak,1} ) plus the conditions
\begin{eqnarray}
&&\f=0\;\;\;\;\mbox{on}\;\;S^-\\
&&{\bf n}\,\cdot\,{\bf A}^2\nabla \f -({\bf v}^2\cdot{\bf n})\,\f=0\;\;
\;\mbox{on}\;\;S^+\;\;.
\label{zadnjass,1}\\
&&\f\;\;\;\;\mbox{is}\;\;L-\mbox{periodic in}\;\;x'=(x_1,\ldots ,x_{n-1}) \;\;.\label{pers,1}
\end{eqnarray}
That suggests the use of method of matched asymptotic expansions (see e.g. \cite{BGM}).\\
We separate the domain in three parts separated by $\Sigma_\ep^+ , \Sigma_\ep^- $ :
\begin{eqnarray*}
&&\Omega^+_\ep = ]\,-L/2 , L/2 \,[^{n-1}\,\times \,]d \ep \log (1 /\ep) , L/2\,[\\
&&\Omega^-_\ep = ]\,-L/2 , L/2 \,[^{n-1}\,\times \,]\,-L/2 ,
 - d \ep \log (1 /\ep)\,[ \\
&&G_\ep = ]\,-L/2 , L/2 \,[^{n-1}\,\times \,]\,-d \ep \log (1 /\ep) , d \ep \log (1 /\ep) \, [ \;\;.
\end{eqnarray*}
Constant $d>0$ is to be determined later in order to minimize the error of approximation.
As suggested, in $\Omega_\ep^\pm $, we approximate $\fe $ by $\f_\varepsilon^0 $ that satisfies the equation (\ref{prvak,1}) and boundary conditions(\ref{rubuv,1}), (\ref{rubuvv,1}) as well as the initial condition (\ref{pocuv,1}).
In $G_\ep $ we look for the asymptotic expansion of $\fe $,  in the form
\begin{eqnarray}
&\fe (x,t) \approx &\f^0_\ep (x,t)+ \ep [ \chi^k_\ep (x/\ep )\;\frac{\partial \f^0_\ep }{\partial x_k } (x,t) +
w_\ep (x/\ep )\;\Phi (t)\; ] +\nonumber\\
&& +\ep^2 [\chi^{k\ell}_\ep (x/\ep )
\frac{\partial^2 \f^0_\ep (x,t)}{\partial x_k\partial x_\ell } +w_\ep^{ij} (x/\ep )\,
\frac{\partial \f^0_\ep (x,t) }{\partial x_i}\,v(x,t)_j+\Phi (t)\,z^k_\ep (x/\ep )\;v(x,t)_k ] +\cdots\;\;. \label{epx,1}
\end{eqnarray}
Here and in the sequel we assume the summation from $1$ to $n$ over the repeating index. 
The function $\f^0_\ep $ mimics the behaviour of $\f $ but has two close jumps in stead of one. In fact, that suggests that more accurate approximation of the real situation would be to have two jumps of the flux; one just above and another just below the array of alveoli. However taking the weak limit smears those two jumps into one.
Namely $\f^0_\ep $ is defined by
\begin{eqnarray}
&&\omega^\ep\;\frac{\partial \f^0_\ep}{\partial t} - \mbox{div}\,({\bf A}^\ep\nabla \f^0_\ep ) +
({\bf v}^\ep\,\cdot\,\nabla\,)\f^0_\ep+\lambda\,\omega^\ep\,\f^0_\ep =0\;\;\;\mbox{in}\;\;
\tilde{\Omega}_\ep^T=(\Omega \backslash (\Sigma^+_\ep \cup \Sigma^-_\ep )\; )\,\times\,]0,T[\nonumber\\
&&\f^0_\ep (x,0)=\f_0 (x)\;\;x\in \tilde{\Omega }_\ep = \Omega \backslash (\Sigma^+_\ep \cup \Sigma^-_\ep ) \nonumber\\
&&\f^0_\ep=0\;\;\;\;\mbox{on}\;\;S^+\\
&&{\bf n}\,\cdot\,{\bf A}^2\nabla \f_\ep^0 -({\bf v}^2\cdot{\bf n})\,\f^0_\ep=0\;\;
\;\mbox{on}\;\;S^-\;\;\nonumber\\
&&[\f^0_\ep ] = 0\;\;\;,\;\;\;\left[{\bf e}_n\,\cdot\,{\bf A}^2\nabla \f^0_\ep
-({\bf v}\cdot{\bf e}_n)\,\f^0_\ep \right] =-\frac{1}{2}\;\Phi\;|\partial P_\ep |\;\;\mbox{on}\;\;\Sigma^\pm_\ep\;\;,\nonumber\\
&&\f^0_\ep\;\;\;\;\mbox{is}\;\;L-\mbox{periodic in}\;\;x'=(x_1,\ldots ,x_{n-1}) \;\;,\nonumber
\end{eqnarray}
with
\begin{eqnarray*}
&&\Sigma^\pm_\ep = ]\,-L/2 , L/2 \,[\,\times \{\pm d\;\ep \log (1 /\ep) \} \;\;\,\,,\\
&&[w](x')=w(x',d\,\ep\,\log(1/\ep )\,+)-w(x',d\,\ep\,\log (1/\ep)\,-)\;\;\;\;x'\in \Sigma^+_\ep\\
&&[w](x')=w(x',-d\,\ep\,\log(1/\ep )\,+)-w(x',-d\,\ep\,\log (1/\ep)\,-)\;\;\;\;x'\in \Sigma^-_\ep
\end{eqnarray*}

 The functions
 $\chi_\ep^k, \chi^{k\ell}_\ep , w^{ij}_\ep $ and $w_\ep $ are the solutions of the auxiliary problems of the stationary diffusion type
 posed in an infinite strip
$${\cal G}_\ep =(\;]-1/2 , 1/2 [^{n-1}\,\times \,{\bf R } \;)\backslash P_\ep\;\;,$$
with
$$P_\ep = M \,\times \,]-\ep^{\beta -1} ,
\ep^{\beta -1}\, [\;\;\;.$$

\begin{picture}(60,40)(-70,130)
\setlength{\unitlength}{1 mm}
\put(100,48){\line(0,1){30}}
\put(112,48){\line(0,1){30}}
\put(103,61.5){\line(1,0){6}}
\put(103,60.5){\line(1,0){6}}
\put(103,60.5){\line(0,1){1}}
\put(109,60.5){\line(0,1){1}}
\put(114,61){\vector(-1,0){5}}

\put(115,61){$\textrm{{\footnotesize $P_\ep$}}$}
\put(95,41){$\textrm{{\footnotesize Figure 2: Strip ${\cal G}_\ep$}}$}\end{picture} \\ 

First two problems read
\begin{eqnarray}
&&-\mbox{div}\,({\bf A}\,\nabla \chi^k_\ep )= 0\;\;\;\mbox{in}\;\;{\cal G}_\ep\nonumber \\
&&{\bf n}\,\cdot\,{\bf A}\,\nabla (\chi^k_\ep +y_k)=0\;\;\;\mbox{on}\;\;\partial P_\ep \label{pom1,1}\\
&&\chi^k_\ep\;\;\;\;\mbox{is}\;\;1-\mbox{periodic in}\;\;y'=(y_1,\ldots ,y_{n-1})  \nonumber\\
&&\lim_{y_n \to \infty }
\nabla \chi^k_\ep = 0\;\;\;.\nonumber
\end{eqnarray}
\begin{eqnarray}
&&-\mbox{div}\,({\bf A}\,\nabla w_\ep )= 0\;\;\;\mbox{in}\;\;{\cal G}_\ep\nonumber \\
&&{\bf n}\,\cdot\,{\bf A}\,\nabla w_\ep=1\;\;\;\mbox{on}\;\;\partial P_\ep \label{pom2,1}\\
&&w_\ep\;\;\;\;\mbox{is}\;\;1-\mbox{periodic in}\;\;y'=(y_1,\ldots ,y_{n-1})   \nonumber\\
&&\lim_{y_n\to\pm\infty} {\bf A}\,\nabla w_\ep (y)=\mp\frac{1}{2}\;|\partial P_\ep |\;{\bf e}_n  \;\;\;.\label{25,1}
\end{eqnarray}

Solvability of problem (\ref{pom1,1}) is classical (see e.g. \cite{Lions} or \cite{Oleinik} ). 
Due to the symmetry of the domain ${\cal G}_\ep $ we obviously have that 
\begin{equation}
\chi^k_\ep (y) = \chi^k_\ep (-y) \;\;.\label{par,1}
\end{equation}
 Furthermore,
there exists a constant $c^k (\ep)$, such that
\begin{equation}
|\chi^k_\ep - c^k (\ep ) |_{H^1 ( |y_n| > s)} \leq C\;e^{-\tau s }\;\;,
\end{equation}
for some $C,\tau > 0 $ .\\
Since $\chi^k_\ep $ is determined up to a constant we may assume in the sequel that $c^k (\ep) =0 $.
\begin{rem}
In general we should have two stabilisation constants $c^k_\pm (\ep) $ at $\pm \infty $. Since (\ref{par,1}) holds those two constants are equal. In case of general $P_\ep $, considered in the first part, this seams not to be the case.
\end{rem}
The problem (\ref{pom2,1}) does not admit a solution with decaying gradient (due to the source term on
$\partial P_\ep $). Therefore we have imposed (\ref{25,1}). To see its behaviour at $\infty $ we need to cut-off that boundary condition first. To do so we take a cut-off function
$$\zeta (y_n)= \left\{
\begin{array}{l}
0\;\;\;\mbox{for}\;\;\;-1/2 <y_n < 1/2\\
1\;\;\;\mbox{for}\;\;\; |y_n |> 1\\
\mbox{smooth otherwise}
\end{array}\right. $$
Now we take
$$\pi (y_n)= -\zeta (y_n)\;\,({\bf A}^{2}_{nn})^{-1}\,| y_n  |\;\frac{1}{2}\;|\partial P_\ep |\;\;\;.$$
The function $ v_\ep (y)=w_\ep (y) - \pi (y_n) $ satisfies the problem

\begin{eqnarray}
&&-\mbox{div}\,({\bf A}\,\nabla v_\ep )= (\mbox{\bf A}_{nn}\pi ')'\;\;\;\mbox{in}\;\;{\cal G}_\ep\nonumber \\
&&{\bf n}\,\cdot\,{\bf A}\,\nabla v_\ep=0\;\;\;\mbox{on}\;\;\partial P_\ep \label{pom3,1}\\
&&v_\ep\;\;\;\;\mbox{is}\;\;1-\mbox{periodic in}\;\;y'=(y_1,\ldots ,y_{n-1})   \nonumber\\
&&\lim_{y_n \to \infty } \nabla v_\ep = 0\;\;\;\nonumber
\end{eqnarray}
and it is obviously pair $v_\ep (y)=v_\ep (-y) $.\\
Since the right-hand side is compactly supported, reasoning as in the case of (\ref{pom1,1}), we conclude that such problem admits  a unique
(up to a constant) solution satisfying
 \begin{equation}
|v_\ep - c(\ep ) |_{H^1 ( |y_n| > s)} \leq C\;e^{-\tau s } \;\;,
\end{equation}
where the constant $c(\ep ) $ can be chosen equal to zero. Therefore the asymptotic behaviour of $w_\ep $, for large $|y_n| $ is
\begin{equation}
 w_\ep (y) \approx -({\bf A}^2_{nn})^{-1}|y_n  |\;\frac{1}{2}\;|\partial P_\ep |+\mbox{exponentially decaying part}\;\;\;.\label{zvj,1}
 \end{equation}
The auxiliary problems for the second corrector are as follows

\begin{eqnarray}
&&-\mbox{div}\,({\bf A}\,\nabla \chi^{\ell m}_\ep )= {\bf A}_{\ell k}\,
\frac{\partial \chi^m_\ep }{\partial y_k}+\frac{\partial }{\partial y_k} ({\bf A}_{k\ell}\,\chi^m_\ep ) \;\;\;\mbox{in}\;\;{\cal G}_\ep\nonumber \\
&&{\bf n}\,\cdot\,{\bf A}\,\nabla \chi^{\ell m} _\ep=0\;\;\;\mbox{on}\;\;\partial P_\ep \label{pomm3,1}\\
&&\chi^{\ell m}_\ep\;\;\;\;\mbox{is}\;\;1-\mbox{periodic in}\;\;y'=(y_1,\ldots ,y_{n-1})  \nonumber\\
&&\lim_{y_n \to \infty }
\nabla \chi^{\ell m}_\ep = 0\;\;\;.\nonumber
\end{eqnarray}

\begin{eqnarray}
&&-\mbox{div}\,({\bf A}\,\nabla w^{ij}_\ep )=\,
\frac{\partial \chi^j_\ep }{\partial y_i} \;\;\;\mbox{in}\;\;{\cal G}_\ep\nonumber \\
&&{\bf n}\,\cdot\,{\bf A}\,\nabla w^{ij} _\ep=0\;\;\;\mbox{on}\;\;\partial P_\ep \label{pom4,1}\\
&&w^{ij}_\ep\;\;\;\;\mbox{is}\;\;1-\mbox{periodic in}\;\;y'=(y_1,\ldots ,y_{n-1})  \nonumber\\
&&\lim_{y_n \to \infty }
\nabla w^{ij}_\ep = 0\;\;\;.\nonumber
\end{eqnarray}

\begin{eqnarray}
&&-\mbox{div}\,({\bf A}\,\nabla z^{k}_\ep )=\,
-\frac{\partial w_\ep }{\partial y_k} \;\;\;\mbox{in}\;\;{\cal G}_\ep\nonumber \\
&&{\bf n}\,\cdot\,{\bf A}\,\nabla z^{k} _\ep=0\;\;\;\mbox{on}\;\;\partial P_\ep \label{pom5,1}\\
&&z^{k}_\ep\;\;\;\;\mbox{is}\;\;1-\mbox{periodic in}\;\;y'=(y_1,\ldots ,y_{n-1})  \nonumber\\
&&\lim_{y_n \to \infty }
\nabla z^{k}_\ep = 0\;\;\;,\;k\neq n\nonumber \\
&&\lim_{y_n\to \infty} ({\bf A}\nabla z^n +\frac{1}{2}(A_{nn}^2)^{-1}\,|\partial P_\ep |\;|y_n|\,)=0  \;\;.\nonumber
\end{eqnarray}

As before we conclude that $\chi^{\ell m}_\ep , w^{ij}_\ep $ and $z^k_\ep\,\;,\; k\neq n $ can be chosen to decay exponentially
towards zero as $ y_n\to \pm \infty $, while
$$z^n_\ep (y)\approx -\frac{1}{4}\;|\partial P_\ep |\;({\bf A}^2_{nn})^{-2}\;|y_n |\,y_n + \mbox{exponentially decaying part}\;\;.$$
\ \\
Now we still have the term $\ep\,w (x/\ep )\;\Phi $ in the inner approximation
that hasn't been matched by our exterior approximation. Taking into account (\ref{zvj,1}) we need to patch our outer approximation with a term
of the following form
$$ d\;\ep\;\log (1/\ep )\;\f^1_\ep \;\;\;.$$
At the same time we will correct the flux jump created by $z^n_\ep $. To do that we define the second corrector $\f^1_\ep $ by

\begin{eqnarray*}
&&\omega^\ep\;\frac{\partial \f^1_\ep}{\partial t} - \mbox{div}\,({\bf A}^\ep \nabla \f^1_\ep ) +
({\bf v}^\ep \,\cdot\,\nabla\,)\f^1_\ep+\lambda\,\omega^\ep \,\f^1_\ep =0\;\;\;\mbox{in}\;\;
\tilde{\Omega}_\ep^T \\
&&\f^1_\ep (x,0)=0\;\;x\in \tilde{\Omega}_\ep\\
&&\f^1_\ep=0\;\;\;\;\mbox{on}\;\;S^+\\
&&{\bf n}\,\cdot\,{\bf A}^2 \nabla \f^1_\ep  -({\bf v}^2\cdot{\bf n})\,\f^1_\ep=0\;\;
\;\mbox{on}\;\;S^-\;\;\\
&&[\f^1_\ep] =\mp\frac{1}{2}\;\Phi\;({\bf A}^2_{nn})^{-1}\;|\partial P_\ep |
\;\;,\;\;\;[{\bf e}_n\,\cdot\,({\bf A}^2 \nabla \f^1_\ep -{\bf v}^2 \,\f^1_\ep )]=\mp\frac{1}{2}\;\Phi\;({\bf A}^2_{nn})^{-1}\;|\partial P_\ep |
\;\;
\;\mbox{on}\;\;\Sigma^\pm_\ep\;\;,\\
&&\f^1_\ep\;\;\;\;\mbox{is}\;\;L-\mbox{periodic in}\;\;x'=(x_1,\ldots ,x_{n-1}) \;\;,
\end{eqnarray*}
For expansion (\ref{epx,1}) we can prove the following error estimate
\begin{tm}
Let $m=\frac{11}{6} $ for $n=3$ and $m<2$ for $n=2$. Let $d\geq 2$. 
There exists a constant $C>0$ independent on $\ep $, such that
\begin{eqnarray}
|\fe - F_\ep  |_{L^2 (0,T;H^1 ({\cal B}_\ep))}\leq C \,(\ep\,\log(1/\ep ))^m\;\;\;.
\end{eqnarray}
where
 $$F_\ep (x,t) = \left\{
\begin{array}{l}
\f^0_\ep (x,t) +d\;\ep\;\log (1/\ep )\;\f^1_\ep (x,t)\;\;\;,\;\;\;\mbox{in}\;\;\;\Omega_\ep^\pm\\
\f^0_\ep (x,t)+ d\;\ep\;\log (1/\ep )\;\f^1_\ep (x,t)+\ep [ \chi^k_\ep (x/\ep )\;\frac{\partial 
}{\partial x_k } (\f^0_\ep +d\,\log (1/\ep)\;\ep\;\f^1_\ep\,) (x,t) + w_\ep (x/\ep )\;\Phi (t)\; ]+\\
+  \ep^2 [\chi^{k\ell}_\ep (x/\ep )
\frac{\partial^2 \f^0_\ep (x,t)}{\partial x_k\partial x_\ell } +w_\ep^{ij} (x/\ep )\,
\frac{\partial \f^0_\ep (x,t) }{\partial x_i}\,v(x,t)_j+\Phi (t)\,z^k_\ep (x/\ep )\;v(x,t)_k ]
 \; \;\;\;,\;\;\;\mbox{in}\;\;\;G_\ep
\end{array}\right.
\;\;.$$

$${\cal B}_\ep = \Omega \backslash (\Sigma^+_\ep \cup \Sigma^-_\ep ) $$
with $\Sigma^\pm_\ep = ]-L/2 ,L/2 [^{n-1}\;\times\;\{\pm\,d\,\ep\log (1/\ep ) \} $. 
Furthermore the same estimate holds in $L^\infty (0,T; L^2 ( \Omega_\ep)) $ norm.\label{tekaj,1}
\end{tm}
{\bf Proof.}
We divide the domain in three parts,
$\Omega^+_\ep , \Omega^-_\ep , G_\ep $.

Let $ R_\ep = \fe - F_\ep $. In $\Omega^\pm_\ep $ we have
$$\omega^2\,\frac{\partial R_\ep }{\partial t} - \mbox{div}({\bf A}^2\nabla R_\ep) + ({\bf v}^2\,\cdot\,\nabla ) R_\ep +
\lambda \omega^2\;R_\ep = 0\;\;.$$
In $G_\ep $ the function $R_\ep $ satisfies
$$\omega^\ep \,\frac{\partial R_\ep }{\partial t} - \mbox{div}({\bf A}^\ep \nabla R_\ep) + ({\bf v}^\ep \,\cdot\,\nabla ) R_\ep +
\lambda \omega^\ep \;R_\ep = E_\ep\;\;,$$
with $|E_\ep |_{L^\infty (G_\ep )}\leq C \,\ep\;\log(1/\ep )\; $ .
Furthermore, on $\Sigma^\pm_\ep $ we have jumps
$$[R^\ep ]= O (\ep^2\log^2 (1/\ep ) ) \;\;\;,\;\;\;[ ({\bf A}^2 \nabla R^\ep  -
R^\ep \,{\bf v}^2 )\,\cdot \,{\bf e}_n \,]=O(\ep^d \,)\;\;\;.$$
Now the result follows by the standard a priori estimate. $\;\; \Box $\\
It should be noticed that the terms in our expansion $F_\ep $ still depend on $\ep $ implicitly.
However it is clear that:
\begin{lm}
\begin{eqnarray*}
&&|\f^0_\ep - \f |_{L^2 (0,T;H^1 (\Omega ) )} \leq C \sqrt{d\;\ep \,\log (1/\ep )} \\
&&|\f^1_\ep |_{L^2 (0,T ; H^1 (\Omega^\pm_\ep )}     \leq C \sqrt{d\;\ep \,\log (1/\ep )} \;\;.
\end{eqnarray*}
\end{lm}
The proof is straightforward and follows by deducing two problems and estimating the remainder.\\
We have the following consequences of theorem \ref{tekaj,1}:
\begin{cor}
Let
$$H_\ep (x) = \left\{
\begin{array}{l}
\f^0_\ep (x) \;\;\;,\;\;\;\mbox{in}\;\;\;\Omega_\ep^\pm\\
\f^0_\ep (x)+ \ep [ \chi^k_\ep (x/\ep )\;\frac{\partial \f^0_\ep }{\partial x_k } (x) + w_\ep (x/\ep )\;\Phi (t)\; ] \;\;\;,\;\;\;\mbox{in}\;\;\;G_\ep^\pm
\end{array}\right.
\;\;.$$
Then
$$|\fe - H_\ep |_{L^2(0,T;H^1 ({\cal B}_\ep))}\leq C (\ep \log(1/\ep )\,)^{3/2}\;\;\;.$$
Furthermore
\begin{eqnarray*}
&&|\varphi_\ep -\varphi |_{L^2 (0,T;H^1 (\Omega_\ep ))}\leq C\sqrt{\ep\;\log(1/\ep)}\\
&&|\varphi_\ep -\varphi^0_\ep |_{L^\infty (0,T;L^2 (\Omega_\ep ))}\leq C\,(\ep\;\log(1/\ep))^{3/2}\;\;.
\end{eqnarray*}
\end{cor}

\ \\
For the auxiliary problems we have:
\begin{lm}
\begin{eqnarray*}
&& |\nabla (\chi^k_\ep - \chi^k )|_{L^2 ({\cal G}_\ep )} \leq C \ep^{(\beta -1) /2 }\\
&& |\nabla (\chi^{mk}_\ep - \chi^{mk} )|_{L^2 ({\cal G}_\ep )} \leq C \ep^{(\beta -1) /2 }\\
&& |\nabla (z^k_\ep - z^k )|_{L^2 ({\cal G}_\ep )} \leq C \ep^{(\beta -1) /2 }\\
&& |\nabla (w_\ep - w )|_{L^2 ({\cal G}_\ep )} \leq C \ep^{(\beta -1) /2 }\\
&& |\nabla (w^{mk}_\ep - w^{mk} )|_{L^2 ({\cal G}_\ep )} \leq C \ep^{(\beta -1) /2 } \,\,,
\end{eqnarray*}
where $ \chi^k ,\chi^{km}, w ,z^k, w^{km} $ are the solutions of corresponding auxiliary
 problems posed on\\
\ \\
\ \\
\ \\
 $${\cal G} = ]-1/2 ,1/2 [\;\times \;{\bf R} \backslash (M\times \{0\} )\;\;.$$
\begin{picture}(60,40)(-70,130)
\setlength{\unitlength}{1 mm}
\put(100,48){\line(0,1){30}}
\put(112,48){\line(0,1){30}}
\put(103,60.9){\line(1,0){6}}
\put(103,61.1){\line(1,0){6}}
\put(114,61){\vector(-1,0){4}}
\put(115,61){$\textrm{{\footnotesize $M\times \{0\}$}}$}
\put(95,41){$\textrm{{\footnotesize Figure 2: Strip ${\cal G}$}}$}
\end{picture} \\ 
\end{lm}
\begin{rem}
It should be noticed that the domain ${\cal G}$ is not locally placed on one side of its boundary. Furthermore, all the problems (\ref{pom1,1}) defining $\chi^k $ for $k\neq n $ have only trivial solutions. On the contrary $\chi^n $ admits a nontrivial solution and the corresponding auxiliary problem now reads

\begin{eqnarray}
&&-\mbox{div}\,({\bf A}\,\nabla \chi^n )= 0\;\;\;\mbox{in}\;\;{\cal G}\nonumber \\
&&{\bf A}_{nk}\frac{\partial (\chi^n +y_n)}{\partial y_k} =0\;\;\;\mbox{on}\;\;M\times \{\pm 0\} \label{pom11,1}\\
&&\chi^n\;\;\;\;\mbox{is}\;\;1-\mbox{periodic in}\;\;y'=(y_1,\ldots ,y_{n-1})  \nonumber\\
&&\lim_{y_n \to \infty }
\nabla \chi^n_\ep = 0\;\;\;.\nonumber
\end{eqnarray}
\begin{eqnarray}
&&-\mbox{div}\,({\bf A}\,\nabla w )= 0\;\;\;\mbox{in}\;\;{\cal G}\nonumber \\
&&\mp{\bf A}_{nk}\,\frac{\partial w}{\partial y_k}=1\;\;\;\mbox{on}\;\;M\times\{\pm 0\} \label{pom22,1}\\
&&w\;\;\;\;\mbox{is}\;\;1-\mbox{periodic in}\;\;y'=(y_1,\ldots ,y_{n-1})   \nonumber\\
&&\lim_{y_n\to\pm\infty} {\bf A}\nabla w = \mp |M|\,{\bf e}_n\;.\nonumber
\end{eqnarray}

\end{rem}
\section{Conclusion}
The expansion (\ref{epx,1}) clearly points out two important terms in the asymptotic behaviour of $\fe $, zero order term $\f^0_\ep $ and first order term $\ep\;w_\ep (x/\ep )\;\Phi $. In a real-life situation, that we are trying to model, the containers are leaking intensively for very short time. During that time $\Phi $ is large and the second order term  $\ep\;w_\ep (x/\ep )\;\Phi $ dominates the behaviour of the solution $\fe $ (despite of $\ep $ multiplying it). Indeed, the typical diffusion coefficient in a low permeable layer (clay) is small, compared to the one in the rest of the domain (limestone). Thus, at the begining of the process, diffusion arround sources is slow. After a short period of time $\Phi $ vanishes and it is only then that the diffusion becomes dominant, i.e. $\f_\ep^0 $ becomes the most important term.

\end{document}